\documentclass[10pt]{article}

\usepackage{fullpage,setspace}

\newcommand{\be}{\begin{equation}} 
\newcommand{\ee}{\end{equation}}
\newcommand{\beq}{\begin{eqnarray}}
\newcommand{\eeq}{\end{eqnarray}}
\newcommand{\nbeq}{\begin{eqnarray*}}
\newcommand{\neeq}{\end{eqnarray*}}

\begin{document}

\title{Characterization of exponential distribution through equidistribution conditions for consecutive maxima}
\date{\empty}

\author{Santanu Chakraborty$^1$ \and George P. Yanev$^1$ \\
}

\maketitle

\footnotetext[1]{Department of Mathematics, The University of Texas - Pan American,
             1201 West \\ \indent \ University Drive, Edinburg, Texas, 78539 USA.}

\begin{abstract}
A characterization of the exponential distribution based on equidistribution conditions for maxima of random samples with consecutive sizes $n-1$ and $n$ for an arbitrary and fixed $n\ge 3$ is proved.
This solves an open problem stated recently in Arnold and Villasenor \cite{AV13}.\\



\end{abstract}

\section{Introduction}
\label{intro}
Characterizations of the exponential distribution are abundant. Comprehensive surveys can be found in Ahsanullah and Hamedani \cite{AH10}, Arnold and Huang \cite{AH95}, and  Johnson, Kotz, and Balakrishnan \cite{JKB94}.  Recently, Arnold and Villasenor \cite{AV13} obtained a series of characterizations based on random sample of size two. They also identified a list of conjectures for possible extensions of their results to bigger samples.  In this work we confirm that one of these conjectures is true for a sample of any fixed size $n\ge 2$. Note that in Yanev and Chakraborty (2013) the case of random sample of size three was considered.

Let $X_1, X_2, \ldots X_n$ for $n\ge 2$ be a random sample from an exponentially distributed parent $X$. It is known that
\be \label{property}
\max\{X_1, X_2, \ldots, X_{n-1}\}+\frac{1}{n}X_n \stackrel{d}{=} \max\{X_1, X_2, \ldots, X_{n}\},
\ee
where $ \stackrel{d}{=}$ denotes equality in distribution. We write $X\sim \exp(\lambda)$ if the probability density function (pdf) of $X$ equals $f_X(x)=\lambda e^{-\lambda x}I(x>0)$.
Our goal is to prove that (\ref{property}), under some regularity assumptions on the cumulative distribution function (cdf) $F$ of $X$, is a sufficient condition for $X$ to be exponential.

\vspace{0.3cm}{\bf Theorem}\ Let $X$ be a non-negative continuous random variable with pdf $f$. If $f$ is analytic in a neighborhood of zero and (\ref{property})
holds true, then $X\sim \exp(\lambda)$ with some $\lambda>0$.

\vspace{0.3cm}Weso{\l}owski and Ahsanullah (2004) and more recently Casta\~{n}o-Martinez et al. (2012) proved characterizations of probability distributions in the context of random translations. Our study differs from these two papers in two ways. First, the condition that the translator random variable has certain known distribution is omitted here. Secondly, a key tool in our proof is a lemma first given in Arnold and Villase\~{n}or (2013). This new direct technique may also be used in obtaining more general results, a possibility which we will explore in the future.

\section{Preliminaries}

 Define for all non-negative integers $m$, $i$, and any real number $x$
\[
H_{m,i}(x):=\sum_{j=0}^m (-1)^j {m \choose j} (x-j)^i.
\]
It is known, (e.g., Ruiz 1996) that for all integers $m\ge 0$ and all real $x$
\be \label{Ruiz}
H_{m,i}(x)=
\left\{
  \begin{array}{ll}
    m! & \mbox{if} \quad i=m; \\
    0 & \mbox{if} \quad 0\le i\le m-1.
  \end{array}
\right.
\ee
\vspace{0.3cm}\noindent The identities given in the next lemma may be of independent interest.

\vspace{0.3cm}{\bf Lemma 1}\ Let $k$, $m$, and $s$ be any positive integers.

\vspace{0.3cm}(i)
\be \label{rec_formula}
H_{m,s}(m+1)=\sum_{j=1}^{s}{s \choose j}H_{m-1,s-j}(m).
\ee

(ii)
\be \label{main_l3}
\sum_{j=0}^k (m+2)^{k-j}H_{m,j}(m+1)=\sum_{j=0}^k {k+1 \choose j+1} H_{m,j}(m+1).
\ee

\vspace{0.3cm}{\bf Proof}. (i)
Indeed,
\nbeq
\lefteqn{\sum_{j=1}^{s}{s \choose j}H_{m-1,s-j}(m)=\sum_{i=0}^{m-1}(-1)^i{m-1 \choose i}\sum_{j=1}^s {s \choose j}(m-i)^{s-j}}\\
    & & =\sum_{i=0}^{m-1}(-1)^i{m-1 \choose i}\left[(m+1-i)^s-(m-i)^s\right]\\
    & & =(m+1)^s-\left[m^s+{m-1 \choose 1}m^s\right]+\left[{m-1 \choose 1}(m-1)^s+{m-1 \choose 2}(m-1)^s\right]\\
   & & \hspace{0.2cm} +\ldots +(-1)^{m-1}\left[{m-1 \choose m-2}2^s+2^s\right]+(-1)^m\\
    & & =(m+1)^s - {m \choose 1}m^s +\ldots +(-1)^{m-1}{m \choose m-1}2^s+(-1)^m\\
    & & =\sum_{j=0}^m (-1)^j{m \choose j} (m+1-j)^s\\
    & & =H_{m,s}(m+1).
    \neeq

(ii) The left-hand side of (\ref{main_l3}) equals
\beq \label{lhs_l3}
\lefteqn{\sum_{j=0}^k (m+2)^{k-j}\sum_{i=0}^m (-1)^i{m \choose i}(m+1-i)^j } \\
    & = &
    \sum_{i=0}^m (-1)^i{m \choose i}(m+2)^k \sum_{j=0}^k \left(\frac{ m+1-i}{m+2}\right)^j \nonumber \\
    & = &
     \sum_{i=0}^m (-1)^i{m \choose i}\frac{1}{i+1}\left[ (m+2)^{k+1}-(m+1-i)^{k+1}\right] \nonumber \\
    & =  &
\sum_{i=0}^m (-1)^i{m+1 \choose i+1}\frac{1}{m+1}\left[ (m+2)^{k+1}-(m+1-i)^{k+1}\right] \nonumber \\
    & = &
-\frac{(m+2)^{k+1}}{m+1}\sum_{r=1}^{m+1}(-1)^r { m+1 \choose  r} + \frac{1}{m+1}\sum_{r=1}^{m+1}(-1)^r { m+1 \choose r}(m+2-r)^{k+1} \nonumber \\
    & = &
-\frac{(m+2)^{k+1}}{m+1}\left[\sum_{r=0}^{m+1}(-1)^r { m+1 \choose  r}- 1\right] \nonumber \\
 & & +\frac{1}{m+1}\left[\sum_{r=0}^{m+1}(-1)^r { m+1 \choose r}(m+2-r)^{k+1}-(m+2)^{k+1}\right] \nonumber \\
    & = &
\frac{1}{m+1}\sum_{r=0}^{m+1}(-1)^r { m+1 \choose r}(m+2-r)^{k+1}. \nonumber
\eeq
For the right-hand side of (\ref{main_l3}) we obtain
\nbeq \label{rhs_l3}
\lefteqn{\sum_{j=0}^k { k+1 \choose j+1}\sum_{i=0}^m (-1)^i{m \choose i}(m+1-i)^j } \\
    & = &
 \sum_{i=0}^m (-1)^i{m \choose i}\sum_{j=0}^k { k+1 \choose j+1}(m+1-i)^j\nonumber \\
    & = &
 \sum_{i=0}^m (-1)^i{m \choose i}\frac{1}{m+1-i}\sum_{j=0}^k { k+1 \choose j+1}(m+1-i)^{j+1}\nonumber \\
    & = &
\frac{1}{m+1}\sum_{i=0}^m (-1)^i{m+1 \choose i}\sum_{r=1}^{k+1} { k+1 \choose r}(m+1-i)^{r}\nonumber \\
    & = &
\frac{1}{m+1}\sum_{i=0}^m (-1)^i{m+1 \choose i}\left[\sum_{r=0}^{k+1} { k+1 \choose r}(m+1-i)^{r}-1\right]\nonumber \\
    & = &
\frac{1}{m+1}\sum_{i=0}^m (-1)^i{m+1 \choose i}(m+2-i)^{k+1}-\frac{1}{m+1} \left[ \sum_{i=0}^{m+1} (-1)^i{m+1 \choose i}-(-1)^{m+1}\right] \nonumber \\
    & = &
\frac{1}{m+1}\sum_{r=0}^{m+1}(-1)^r { m+1 \choose r}(m+2-r)^{k+1}, \nonumber
\neeq
which equals (\ref{lhs_l3}). The proof of the lemma is complete.

\vspace{0.5cm}{\bf Lemma 2}\ Assume $F(0)=0$. Let $m$ and $d$ be integers, such that $m>0$ and $d\ge -m$.
Then the $m+d$th derivative of $G_m(x):=F^m(x)f(x)$ at 0 is given by
\be \label{main}
G_m^{(m+d)}(0)=
\left\{
  \begin{array}{ll}
  \left(\frac{f'(0)}{f(0)}\right)^{d}f^{m+1}(0)H_{m,m+d}(m+1) & \mbox{if} \quad d\ge 0;
 \\
    0 & \mbox{if} \quad -m\le d<0,
  \end{array}
\right.
\ee
where for $d>0$ we assume additionally that
\be \label{der}
f^{(k)}(0)= \left(\frac{f'(0)}{f(0)}\right)^{k-1}f'(0) \qquad k=0,1,\ldots
\ee

\vspace{0.3cm}{\bf Proof}. We will consider separately the cases: (i) $-m\le d<0$; (ii) $d=0$, and (iii) $d>0$.

(i) If $-m\le d<0$, then $G^{(m+d)}_m(0)=0$ because all the terms in the expansion of
$G^{(m+d)}_m(0)$ has a factor $F(0)=0$.

(ii) Let $d=0$. We shall prove (\ref{main}) by induction on $m$. One can verify directly (\ref{main}) for $m=1$.
Assuming (\ref{main}) for $m=k$, we shall prove it for $m=k+1$.
Since $G_{k+1}(x)=F(x)G_{k}(x)$, we have
\nbeq
G^{(k+1)}_{k+1}(0) & = & \sum_{j=0}^{k+1} {k+1 \choose j}F^{(j)}(0) G^{(k+1-j)}_k(0)\\
    & = & F(0)G^{(k+1)}_k(0)+(k+1)F'(0)G^{(k)}_k(0)\\
    & = & (k+1)!f^{k+2}(0),
\neeq
which completes the proof in this case.

(iii) We turn to the proof of (\ref{main}) for $d>0$ and any positive integer $m$.
 For simplicity, in this part of the lemma's proof, we will write $f:=f(0)$ and $f^{(j)}:=f^{(j)}(0)$.

(a) Let $m=1$. If $d=1$, then we have $G_1^{(2)}(0)=3f'f$ and
\[
f'fH_{1,2}(2)=f'f\sum_{j=0}^1(-1)^j{1 \choose j}(2-j)^2=3f'f.
\]
Thus, (\ref{main}) is true for $d=1$. Next, assuming that (\ref{main}) hold true for $G_1^{(k)}(0)$, that is $d=k-1\ge 1$, we shall prove it for $G_1^{(1+k)}(0)$, that is $d=k$.
Indeed, using the assumption (\ref{der}) we obtain
\nbeq
G_1^{(1+k)}(0) & = & \sum_{j=0}^{k+1} {k+1 \choose j} F^{(j)}f^{(k+1-j)}\\
    & = & \sum_{j=1}^{k+1} {k+1 \choose j} f^{(j-1)}f^{(k+1-j)}\\
    & = & \sum_{j=1}^{k+1} {k+1 \choose j} \left(\frac{f'}{f}\right)^{j-2}f'\left(\frac{f'}{f}\right)^{k-j}f'\\
    & = & \left(\frac{f'}{f}\right)^{k-2}(f')^2\sum_{j=1}^{k+1}{k+1 \choose j} \\
    & = & \left(\frac{f'}{f}\right)^{k-2}(f')^2(2^{k+1}-1).
    \neeq
Using (\ref{der}) again and the definition of $H_{1,1+k}(2)$, for the right-hand side of (\ref{main}) we have
\nbeq
\left(\frac{f'}{f}\right)^{k}f^{2}H_{1,1+k}(2) & = & \left(\frac{f'}{f}\right)^{k-2}(f')^2 \sum_{j=0}^1 (-1)^j {1 \choose j}(1-j)^{k+1}\\
    & = & \left(\frac{f'}{f}\right)^{k-2}(f')^2 (2^{k+1}-1).
    \neeq
This completes the proof of the induction step for the case (a) $m=1$ and any $d>0$.

(b) Assuming (\ref{main}) for $m=1, 2, \ldots k$ and any $d>0$ we shall prove it for $m=k+1$ and any $d>0$.
Using (\ref{der}) and the induction assumption, we obtain
\nbeq
G_{k+1}^{(k+1+d)}(0)
    &  = &
\sum_{j=1}^{k+1+d} {k+1+d \choose j} f^{(j-1)}G_{k}^{(k+1+d-j)}(0)\\
 & = & \sum_{j=1}^{k+1+d}{k+1+d \choose j}\left(\frac{f'}{f}\right)^{j-2}f'G_{k}^{(k+1+d-j)}(0)\\
    & = &
    \sum_{j=1}^{k+1+d}{k+1+d \choose j}\left(\frac{f'}{f}\right)^{j-2}f'\left(\frac{f'}{f}\right)^{1+d-j}
    f^{k+1}H_{k,k+1+d-j}(m)\\
    &  = &
    \left(\frac{f'}{f}\right)^{d}f^{k+2} \sum_{j=1}^{k+1+d}{k+1+d \choose j}H_{k,k+1+d-j}(k+1)\\
 &   = &
\left(\frac{f'}{f}\right)^{d}f^{k+2} H_{k+1,k+1+d}(k+2),
    \neeq
where the last equality follows from (\ref{rec_formula}) with $s=k+1+d$. This proves the induction step (b).
It follows from (a) and (b) that the statement in (iii) is true. The lemma is proven.

\vspace{0.3cm}Next lemma, due to Arnold and Villase\~{n}or \cite{AV13}, plays a crucial role in the proof of the theorem.  For completeness of the exposition, we provide the proof too.

\vspace{0.5cm}{\bf Lemma 3}\ If $F(0)=0$, the pdf $f$ is analytic in a neighborhood of 0, and
\be \label{lemma}
f^{(k)}(0)=\left(\frac{f'(0)}{f(0)}\right)^{k-1}f'(0), \qquad k=0,1,\ldots,
\ee
then $X\sim \exp\{\lambda\}$ for some  $\lambda>0$.

\vspace{0.3cm}{\bf Proof.}\ For the Maclaurin series of $f(x)$, we have for $x>0$
\beq \label{Mcl}
f(x) & = & \sum_{k=0}^\infty \frac{f^{(k)}(0)}{k!}x^k \\
    & = & f(0) + \sum_{k=1}^\infty \left(\frac{f'(0)}{f(0)}\right)^{k-1}f'(0)\frac{x^k}{k!} \nonumber \\
    & = &
    f(0)\exp\left\{\frac{f'(0)}{f(0)}x\right\}. \nonumber
\eeq
Since $f(x)$ is a pdf, we have $f'(0)/f(0)<0$. Denoting $\lambda=-f'(0)/f(0)>0$ and setting the integral of (\ref{Mcl}) from 0 to $\infty$ to be 1, we obtain $\lambda=f(0)$. Therefore, $f(x)=\lambda e^{-\lambda x}I(x>0)$, i.e., $X\sim \exp\{\lambda\}$.

\section{Proof of the theorem}

Equation (\ref{property}) can be written in terms of the probability densities of the random variables involved as follows
\[
\int_0^x f_{X_n/n}(y)f_{\max\{X_1,\ldots,X_{n-1}\}}(x-y)\, dy=nf(x)F^{n-1}(x).
\]
This is equivalent to
\[
\int_0^x nf(ny)(n-1)F^{n-2}(x-y)f(x-y)\, dy=nf(x)\int_0^x(n-1)F^{n-2}(y)f(y)\, dy,
\]
which simplifies to
\be \label{eqn8}
\int_0^x f(ny)G_{n-2}(x-y)\, dy=f(x)\int_0^xG_{n-2}(y)\, dy,
\ee
where $G_m(x):=F^m(x)f(x)$, as before.
In view of Lemma 3, to complete the proof of the theorem it  suffices to show
\be  \label{eqn16}
f^{(k)}(0)=\left(\frac{f'(0)}{f(0)}\right)^{k-1}f'(0), \qquad k=0,1,\ldots
\ee
The equation (\ref{eqn16}) trivially holds for $k=0$ and $k=1$. Let $r$ be an arbitrary positive integer. Assuming (\ref{eqn16}) for all $1\le k\le r$, we shall prove it for $k=r+1$.
Differentiating the left-hand side of (\ref{eqn8}) $(n+r)$ times, we obtain
\beq  \label{eqn9}
\lefteqn{\frac{d^{n+r}}{dx^{n+r}}\int_0^xf(ny)G_{n-2}(x-y)\, dy}\\
    & = & \sum_{i=0}^{n+r-1} n^{n+r-1-i}f^{(n+r-1-i)}(nx)G_{n-2}^{(i)}(0)+\int_0^x f(ny)G_{n-2}^{(n+r)}(x-y)\, dy. \nonumber
\eeq
Applying Leibnitz product rule of differentiation to the right-hand side, we have
\beq \label{eqn10}
\lefteqn{\frac{d^{n+r}}{dx^{n+r}}\left[f(x)\int_0^xG_{n-2}(y)\, dy\right]}\\
    & = &
\sum_{i=0}^{n+r-1}{n+r \choose i+1}f^{(n+r-1-i)}(x)G_{n-2}^{(i)}(x)+f^{(n+r)}(x)\int_0^x G_{n-2}^{(n+r)}(y)\,dy \nonumber
\eeq
Setting $x=0$ in (\ref{eqn9}) and (\ref{eqn10}), we see that (\ref{eqn8}) implies
\[
\sum_{i=0}^{n+r-1} a_if^{(n+r-1-i)}(0)G_{n-2}^{(i)}(0)=\sum_{i=0}^{n+r-1}b_if^{(n+r-1-i)}(0)G_{n-2}^{(i)}(0),
\]
where
\[
a_i:=n^{n+r-1-i}\qquad \mbox{and}\qquad b_i:={n+r \choose i+1}.
\]
Using the induction assumption and the formulas for $G_{n-2}^{(i)}(0)$ in Lemma~2, we write the last equation as
\[
\hspace{-0.3cm}\sum_{i=n-2}^{n+r-1}\! a_if^{(n+r-1-i)}(0)
\!\left(\frac{f'(0)}{f(0)}\right)^{i}\!H_{n-2,i}(n-1)\!
= \! \sum_{i=n-2}^{n+r-1}\! b_if^{(n+r-1-i)}(0)
\! \left(\frac{f'(0)}{f(0)}\right)^{i}\!H_{n-2,i}(n-1).
\]
Separating the summands with $i=n-2$ from the rest of the sums, we have
\nbeq
\lefteqn{a_{n-2}f^{(r+1)}(0)H_{n-2,n-2}(n-1)+\sum_{i=n-1}^{n+r-1} a_{i}f^{(n+r-1-i)}(0)
\left(\frac{f'(0)}{f(0)}\right)^{i}H_{n-2,i}(n-1)}\\
& = &
\hspace{-0.3cm}b_{n-2}f^{(r+1)}(0)H_{n-2,n-2}(n-1)+\sum_{i=n-1}^{n+r-1} b_{i}f^{(n+r-1-i)}(0)
\left(\frac{f'(0)}{f(0)}\right)^{i}H_{n-2,i}(n-1).
\neeq
Set $c_i:=a_i-b_i$ and note that for $n-1\le i\le n+r-1$ we have $0\le n+r-1-i\le r$. Thus applying the induction assumption to
$f^{(n+r-1-i)}(0)$, we obtain
\[
c_{n-2}f^{(r+1)}(0)H_{n-2,n-2}(n-1)
    +
    \left(\frac{f'(0)}{f(0)}\right)^{r}f'(0)\sum_{i=n-1}^{n+r-1}c_{i}H_{n-2,i}(n-1)=0.
\]
To establish (\ref{eqn16}) for $k=r+1$ we need to prove
\[
c_{n-2}H_{n-2,n-2}(n-1)
     + \sum_{j=n-1}^{n+r-1}c_{i}H_{n-2,i}(n-1)=0
\]
or, taking into account that $H_{n-2,i}(n-1)=0$ for $0\le i\le n-3$,
\be \label{last}
\sum_{j=0}^{n+r-1}n^{n+r-1-i}H_{n-2,i}(n-1)=
     \sum_{j=0}^{n+r-1}{n+r \choose i+1}H_{n-2,i}(n-1).
\ee
The equation (\ref{last}) follows from Lemma~1(ii). This completes the induction step and thus proves (\ref{eqn16}). Referring to (\ref{eqn16}) and Lemma~3 we complete the proof of the theorem.


\end{document}